\documentclass[11pt]{article}

\usepackage{latexsym}

\usepackage{amsmath}
\usepackage{amssymb}

\usepackage{graphicx}

\usepackage{color,colortbl}

\def\be{\begin{equation}}
\def\ee{\end{equation}}
\def\ba{\begin{eqnarray*}}
\def\ea{\end{eqnarray*}}
\def\baa{\begin{eqnarray}}
\def\eaa{\end{eqnarray}}

\def\b{\beta}
\def\cl{\centerline}
\def\md{\mbox{\rm mid\,}}
\def\rd{\mbox{\rm rad\,}}
\def\R{{\Bbb R}}
\def\eps{\varepsilon}

\title{\Large\bf A note on non-inclusion ``interval" root solvers}

\author{{\bf Miodrag S. Petkovi\'c$^{\,1,}$, Ljiljana D. Petkovi\'c$^{\,2,}$}\footnote{Corresponding author}
   \\[2mm]
   $^1$\small\it  Faculty of Electronic Engineering,
University of Ni\v s,\\[-1mm]
\small\it A. Medvedeva 14, 18000 Ni\v s, Serbia\\
$^2$\small\it  Faculty of Mechanical Engineering,
University of Ni\v s, \\[-1mm] \small \it A. Medvedeva 14, 18000 Ni\v s, Serbia}

\date{}

\begin{document}

\maketitle

\begin{abstract}
In this note we give comments on recent interval extensions of the King-Ostrowski multipoint methods for solving nonlinear equations published in at least five scientific journals offering faulty ``self-validated" methods. We show that they do not  possess essential feature of genuine interval methods -- inclusion property.
\\[1mm]
\end{abstract}

\renewcommand{\thefootnote}{}%
\footnote{{\it E-mail address:}  ljiljana.petkovic@masfak.ni.ac.rs
 }

\section{Introduction}

In this note we wish to draw attention to erroneous treatment of interval methods  for the inclusion of an isolated zero of a given real function, published in the recent papers \cite{lotfi-1}--\cite{eft-numa} and maybe in some other journals or proceedings, unavailable to the author of this note. Using a simple derivation of relations necessary for the development of genuine inclusion methods (also called interval or self-validated methods), examples and counterexamples, we point to the key mistakes made in the cited papers. Beside these comments,  we also wish to prevent possible construction of new  methods that suffer from the lack of  the most important feature of interval methods -- inclusion property. The mentioned mistakes are mainly made due to the careless and incorrect application of properties of interval arithmetics and the mean value theorem.
 The basic operations and properties of real interval arithmetic can be found in the books \cite{AH}--\cite{mur3}  and many papers cited therein so that we do not list them in this note.

\smallskip

As well known, the main advantage of  interval methods for the enclosure of an isolated zero $x^*$ of a given function is the ability to produce a sequence of nested intervals
$X_0\subset X_1\subset X_2\subset \cdots$ such that $x^*\in X_k$ for all $k=0,1,\ldots\;.$ Methods start from an initial  interval $X_0\owns x^*$ and terminate when the given stoping criterion is satisfied.

\smallskip

Let $X=[\underline{x},\overline{x}]\;(\underline{x},\overline{x}\in \R,\;\underline{x}\le\overline{x}])$ be a closed real interval. The set of all closed real intervals will be denoted by $I(\R)$ and the members of this set by capital letters $A,B,C,\ldots\ .$ For a real function $f,$ define its {\it interval extension}  $f(X)=\{f(x)\;|\;x\in X\}$  that satisfies the conditions
\baa
&&f([x,x])=f(x)\quad  (\mbox{\rm restriction}),\nonumber\\
&&f(x)\in f(X)\quad (\mbox{\rm inclusion property})\label{1}.
\eaa

The goal of the papers \cite{lotfi-1}--\cite{eft-numa} was to extend Ostrowski's and King's two-point methods (applicable in $\R$) to the set $I(\R).$ Both methods can be expressed by the unique iterative formula
\be
\left\{\begin{array}{l}
$$y_k=x_k-\dfrac{f(x_k)}{f'(x_k)},\\[10pt]
x_{k+1}=y_k-\dfrac{f(x_k)+\b f(y_k)}{f(x_k)+(\b-2)f(y_k)} \cdot \dfrac{f(y_k)}{f'(x_k)},
\end{array}\quad (k=0,1,2,\ldots),
\right. \label{2}
\ee
where $\b$ is a real parameter. Formula (\ref{2}) defines two-point King's method \cite{king} of order four.  Ostrowski's method
\cite{ostrov} is obtained from (\ref{2}) as a special case for $\b=0.$

\smallskip

Let $X=[\underline{x},\overline{x}],$ then we define the {\it midpoint} and {\it semi-width} of $X$ as
$$\md(X)=\frac{\underline{x}+\overline{x}}{2}, \ \  \rd(X)=   \frac{\overline{x}-\underline{x}}{2},$$ respectively.
For two intervals $X_k$ and $Y_k$ let us introduce the abbreviation
$$
t_k=\frac{f\bigl(\md(Y_k)\bigr)}{f\bigl(\md(X_k)\bigr)}.
$$

\smallskip
The following  two-point method was constructed in \cite{lotfi-1}, \cite{ferara}, \cite{eft-numa} as the interval extension of (\ref{2})

\medskip
\hskip3.3cm  Given $X_0\owns x^*,$ for $k=0,1,2,\ldots$ calculate
\vspace{-1mm}
\be
\left\{\begin{array}{l}
Y_k=N(X_k)\cap X_k,\\
K(X_k,Y_k)=\md(Y_k)-\dfrac{1+\b t_k}{1+(\b-2)t_k}
\cdot \dfrac{f\bigl(\md(Y_k)\bigr)}{F'(X_k)},\\[7pt]
X_{k+1}=K(X_k,Y_k)\cap X_k,
\end{array}
\right.\label{3}
\ee
referred as {\it Interval Ostrowski's method} for $\b=0,$  and {\it Interval King's method} for $\b\in \R$, see \cite{fazi}, \cite{eft-AMC}. By the same approach the authors constructed  the three-point method

\medskip
\hskip2.9cm  Given $X_0\owns x^*,$ for $k=0,1,2,\ldots$ calculate
\vspace{-1mm}
\be
\left\{\begin{array}{l}
Y_k=N(X_k)\cap X_k,\\
K(X_k,Y_k)=\md(Y_k)-\dfrac{1+\b t_k}{1+(\b-2)t_k}
\cdot \frac{f\bigl(\md(Y_k)\bigr)}{F'(X_k)},\\[7pt]
Z(X_k)=K(X_k,Y_k)\cap X_k,  \\
M(X_k,Y_k,Z_k)=\md(Z_k)-\dfrac{1+\b t_k}{1+(\b-2)t_k}
\cdot \dfrac{f\bigl(\md(Z_k)\bigr)}{F'(X_k)},\\
X_{k+1}=M(X_k,Y_k,Z_k)\cap X_k,\\
\end{array}
\right.\label{4}
\ee
referred as {\it Interval modified Ostrowski method} for $\b=0$ in  \cite{ferara}, \cite{eft-numa}.

\section{Construction of self-validated methods}

Mistakes given in the cited papers \cite{lotfi-1}--\cite{eft-numa} can be  divided into two groups:

\smallskip
(I) The loss of inclusion property due to the wrong interpretation of the mean value theorem.

\smallskip
(II) The wrong estimation of the convergence order of the proposed interval methods as a consequence of incorrect application of properties of interval arithmetics and negligent estimation procedure.

\smallskip

First, we analyze the  group (I) of mistakes that are of decisive importance.
 As in the papers \cite{lotfi-1}--\cite{eft-numa},  in this note we also consider only real-valued functions $f$ differentiable on a specific interval $X=[a,b]\in I(\R)$ such that $f(a)f(b)<0.$ This means that there is at least one zero $x^*$ in $X.$ Furthermore, assume that $f$ is monotonically increasing (decreasing) on $X$, that is, $f'(x)>0$ ($f'(x)<0$) for all $x\in X$.
 According to the previous conditions it follows that there exists only one  zero $x^*$ in $X$ and $0\notin f'(X).$

\smallskip

  For two numbers $a,b\in \R$ define the interval $I(a,b):=[\min\{a,b\},\max\{a,b\}]$ and let $I(x,x^*)\subseteq X\ (x,x^*\in X).$ According to the mean value theorem (shorter MVT) we have
\be
f(x)=f(x^*)+(x-x^*)f'(\xi),\ \ \xi\in I(x,x^*),\label{5}
\ee
or (since $f(x^*)=0$)
\be
f'(\xi)=\frac{f(x)}{x-x^*}. \label{6}
\ee
Simple geometrical interpretation of MVT shows that given any chord of a smooth curve, we can find a point (denoted here with $\xi$) lying between the end-points of the chord such that the tangent at that point is parallel to the chord.
From (\ref{6}), by using inclusion property (\ref{1}), we obtain
\be
x^*=x-\frac{f(x)}{f'(\xi)}\in x-\frac{f(x)}{f'(I(x,x^*))}\subseteq  x-\frac{f(x)}{F'(X)}=:N(X),\label{7}
\ee
where $N$ is the Newton operator.
$F'(X)$ is an interval extension of  $f'(I(x,x^*))$, that is, $F'(X)\supseteq  f'(I(x,x^*)).$ Clearly, we are forced to use $F'(X)$ since $f'(I(x,x^*))$ is unknown interval since $x^*$ is unknown. Note that $F'(X)=I\bigl(f'(\underline{x}),f'(\overline{x})\bigr)$ since $f$ is monotone on $X=[\underline{x},\overline{x}].$

\smallskip

Taking $x=\md(X)$ in (\ref{7}) R. E. Moore constructed in \cite{mur1} (see, also, \cite{AH}--\cite{mur3}) the following inclusion method of Newton's type:

\medskip
\hspace{4.5cm}Given $X_0\owns  x^*,$ for $k=0,1,2,\ldots$ calculate
\vspace{-2mm}
\baa
\left\{\begin{array}{l}
N(X_k)=\md(X_k)-\dfrac{f\bigl(\md(X_k)\bigr)}{F'(X_k)}\,\\[7pt]
X_{k+1}=N(X_k)\cap X_k.
\end{array}
\right. \label{8}
\eaa

The Moore-Newton interval method (\ref{8}) is presented here only to demonstrate
essential characteristics of inclusion methods.
Geometrical interpretation of the  interval method (\ref{8}) is displayed on Fig. 1 (shaded area and the interval $Y$), where the tangents at the endpoints are translated to the midpoint $m(X).$ For clear demonstration, the simple function $f(x)=x^2-4$ and the initial interval $X=[1,4]\owns x^*=2$ have been taken.

\cl{\includegraphics[height=5.7cm]{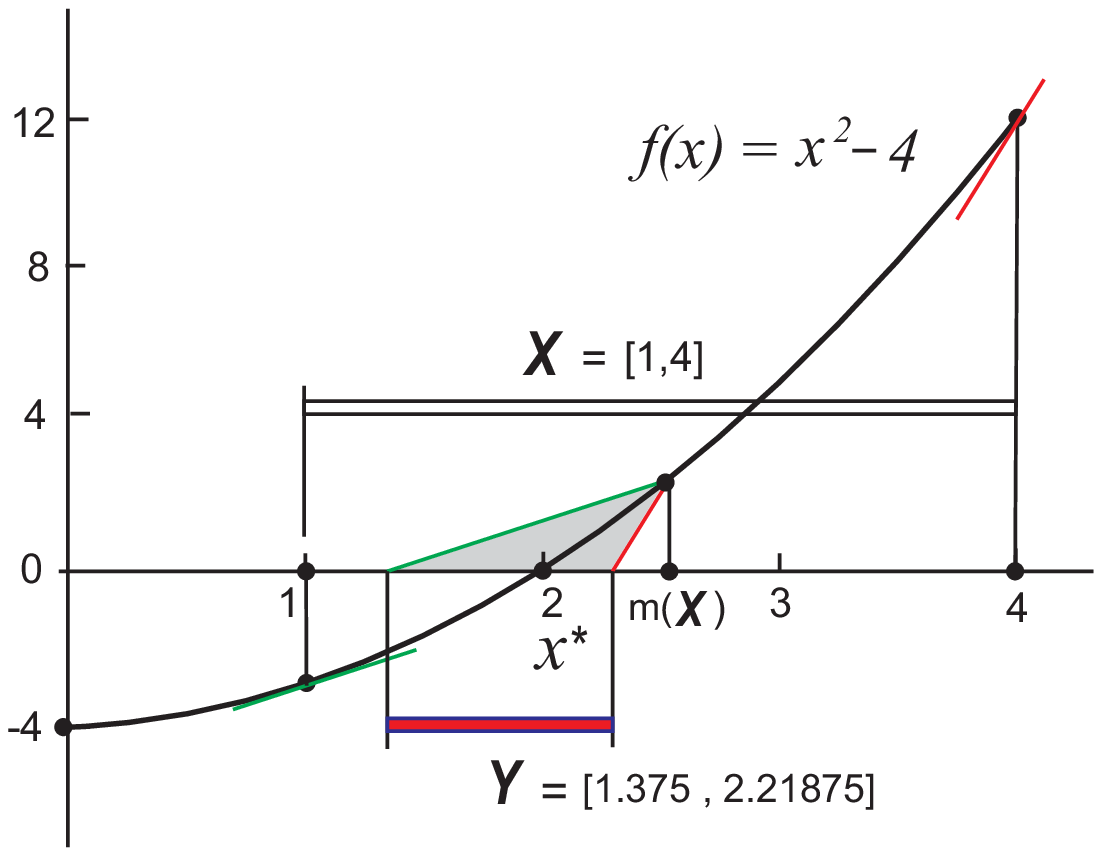}}

\cl{\small Fig. 1 The Moore-Newton interval method: geometrical interpretation}

\section{Lost inclusion property and incorrect convergence order}

The key mistake in the construction of the methods (\ref{3}) and (\ref{4}), which causes the lack of the inclusion property, is the use a ``variation" of
MVT (\ref{5}) (given by (\ref{9})) instead of MVT (\ref{5}).
In fact, the authors of the cited papers   started from the relation
 \be
f(x)=f(x^*)+(x-x^*)c f'(\xi),\ \ \xi\in I(x,x^*),\label{9}
\ee
where $c\neq 1$ is a real number that depends on the values of the function $f$ evaluated at two points belonging to two successive intervals from the sequence $\{X^{(k)}\}.$ Since $f(x^*)=0,$ from (\ref{9}) we obtain
\be
f(\xi')=\frac{f(x)}{c(x-x^*)}. \label{10}
\ee

The three-point method (\ref{4}) has a similar form as the two-point method (\ref{3}) so that we will analyze only the method (\ref{3}). For simplicity, we will omit the iteration index $k.$ Without loss of generality, we will deal with arbitrary points $x\in X=[\underline{x},\overline{x}]$ and $y\in Y=(\underline{y},\overline{y}),$ including $x=\md(X)$ and $y=\md(Y)$ (standing in (\ref{3})). The value  $c$ appearing in (\ref{3}) and  (\ref{4}) is given by
$$
c=\frac{1+(\b-2) t}{1+\b t},\ \ t=\frac{f(y)}{f(x)}\ \ \ (\b\in \R).
$$

\cl{\includegraphics[height=7.5cm]{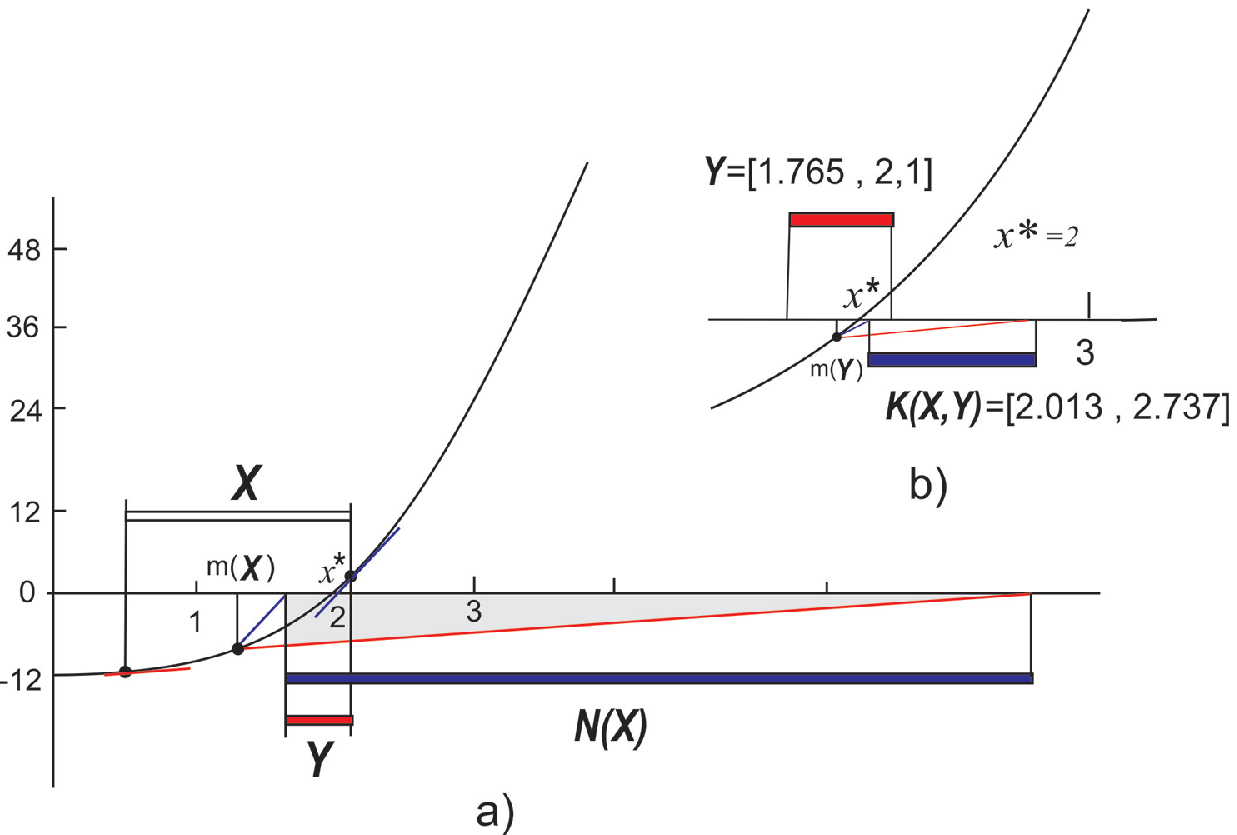}}

\cl{\small Figure 2 The failure of Ostrowski-like ``interval" method}

\medskip

Let us consider the second step of (\ref{3}) (suppressing the index $k$) in the form
\be
K(X,Y)=x-\frac{1+\b t}{1+(\b-2)t}\cdot \frac{f(y)}{f'(I(\underline{x},\overline{x}))}=x- \frac{f(y)}{c f'(I(\underline{x},\overline{x}))}
=x- \frac{f(y)}{I\bigl(c f'(\underline{x}\bigr),c f'(\overline{x})\bigr)}. \label{11}
\ee
Since $t=f(y)/f(x)\ne 0$ there follows $c\ne 1.$ Hence we observe that the slopes $cf'(\underline{x})$ and $c f'(\overline{x})$  {\it do not coincide with the tangents at the end-points} $\underline{x}$ and $\overline{x}$, used in (\ref{8}). The new interval $K(X,Y)$  is now {\it shifted}  and {\it it does not have to contain the zero} $x^*$, as displayed in Fig. 2b. In essence, MVT (\ref{5}),
which {\it preserves} the inclusion property, differs from the relation (\ref{9}) which leads to the described unfavorable situation. As a consequence, the {\it inclusion property of the method {\rm (\ref{3})} is not guaranteed.}

\smallskip
  Geometrical interpretation (taking the function $f(x)=x^3+x^2-12$ and the interval $[0.5,2.1]$ for demonstration) is displayed in Fig. 2. Observe that the new interval $K(X,Y),$ constructed by green and red  lines, does not contain the zero $x^*=2.$  Numerical data are given in Example 1. For the presented example the method (\ref{3}) fails. The same would happen with the method (\ref{4}).

\medskip

{\bf Remark 1.} According to the results presented in the papers \cite{lotfi-1}--\cite{eft-numa},  we conclude that the proofs of inclusion property are either wrong (see \cite{lotfi-1}, \cite{eft-numa}, \cite{ferara}) or  omitted (see \cite{fazi}, \cite{eft-AMC}). The wrong proofs are the consequence of careless and incorrect manipulations with the basic arithmetic operations, incorrect relations and inequalities,  as well as the use of ``variation" of MVT. These mistakes are elementary so that we skip details.

\medskip

{\bf Remark 2.} The basic (wrong) idea presented in \cite{eft-AMC} (published 2015) is  the same as  in \cite{fazi} and \cite{ferara}. Furthermore, a similar (wrong) idea exposed in \cite{eft-numa} (2015)  appeared previously in \cite{ferara} (2013). It is hard to believe that   the same methods are discovered independently since the latter authors cited the former authors.  However, in this note we do not intend to discuss ethical principles. By the way, such kind of job is not necessary since the results of all papers \cite{lotfi-1}--\cite{eft-numa} are wrong.

\medskip

{\bf Remark 3.} To provide inclusion intervals in each iteration, in \cite{lotfi-1} and \cite{ferara}, the authors need to check the conditions (i) $K(X_k,Y_k)\subset X_k\ (\owns x^*)$ for (\ref{3}) and (ii) $\{K(X_k,Y_k)\subset X_k\ (\owns x^*)\ \wedge\  M(X_k,Y_k,Z_k)\subset X_k\ (\owns x^*)\}$ for (\ref{4}) {\it in each iteration}. This is unacceptable iteration model, not only very expensive but also  contrary to standard iterative processes.
What will happen if (i) and (ii) do not hold for some $k$?
\medskip

{\bf Remark 4.} As numerical examples presented in \cite{lotfi-1}--\cite{eft-AMC} have shown, it is possible to run the methods (3) and (4) in such a way that they produce intervals containing the sought zero. Such situations exist in the case when the factor $c$ in (\ref{3}) and (\ref{4}) is close to 1, which will happen when $\rd(X)$ is sufficiently small so that the value $t=f(y)/f(x)=O(\rd(X))$ is also small in magnitude.  Then the dislocation of the shifted interval $cF'(X)=I(c f'(\underline{x}),c f'(\overline{x})$ may be small enough so that the intervals $K(X,Y)$ and $M(X,Y,Z)$ enclose the zero $x^*$. However, since $c\ne 1$ (although $c$ can be very close to 1), using multiple-precision arithmetic in order to obtain inclusion intervals of very small widths it is possible that the zero $x^*$ comes out from some interval in later iterations, see Example 3. Such unfavorable case could not have  happened if a proper interval method had been applied (e.g., the Moore-Newton method (\ref{8})). Simply, as mentioned in Remark 1, the methods (3) and (4) cannot guarantee the enclosure of a zero in each iteration.

\medskip
Among many counterexamples we present two of them which, together with the above discussion and Fig. 2,  confirm that the methods published in  \cite{lotfi-1}--\cite{eft-numa} do not guarantee the enclosure of the sought zero by produced intervals.

\medskip

{\bf Example 1.} We take $\b=0$ in (\ref{3}), that is, we wish to iterate ``interval" Ostrowski-like method (\ref{3}).  Let us find the enclosure for the zero $x^*=2$ of the function $f(x)=x^3+x^2-12$ on the interval $X_0=[0.5,2.1].$ Since $f'(x)=3x^2+2x>0$ for all $x\in X_0,$ $f$ is monotonically increasing on this interval and $0\notin F'(X_0).$ This means that the division by zero-interval does not appear in  (\ref{3}) and (\ref{4}). First, we find
$$
N(X_0)=\md(X_0)-\dfrac{f(\md(X_0))}{F'(X_0)}=1.3-\frac{f(1.3)}{[f'(0.5),f'(2.1)]}
=1.3-\frac{-8.113}
{[1.75,17.43]}=[1.76546,5.936].
$$
Hence $x^*\in N(X_0)$ (the inclusion holds, see Fig. 2a) and the process can be continued. First, we find
$$
Y_0=N(X_0)\cap X_0=[1.76546,5.936]\cap[0.5,2.1]=[1.76546,2.1]
$$
and calculate in interval arithmetic
$$
K(X_0,Y_0)=\md(Y_0)-\frac{1}{1-2f(\md(Y_0))/f(\md(X_0))}\cdot \frac{f(\md(Y_0))}{F'(X_0)}
=[2.01348,2.73702]
$$
(see Fig. 2b). We observe that $x^*=2\notin K(X_0,Y_0))$ and the Ostrowski-like method  fails  at the beginning of iterative process.

\medskip

{\bf Example 2.} King-like method (\ref{4})  breaks down for the function $f(x)=x^3-8,$ the initial interval $X_0=[1.5,2.3]$ and  $\b=5.$ One obtains $K(X_0,Y_0)=[2.024393,2.029699]\not\owns x^*=2.$

\medskip

{\bf Example 3.} To examine the extent of validity of the King-like method (\ref{3}),  we have tested the monic polynomials of the 7th degree of the form
$$
P(x)=(x-1)(x^6+a_5 x^5+a_4 x^4+a_3 x^3+a_2 x^2+a_1 x+a_0),
$$ where the coefficients $a_0,\ldots,a_5$ have been chosen randomly from the interval (0,1) using the statement {\tt Random[\;]} in computer algebra system  {\it Mathematica}.   The initial interval $X_0=[0.6,3]$ containing the zero $x^*=1$ has been chosen. We have tested 100 random polynomials with the values of  the parameter $\b=-2+0.5m\ (m=0,1,\ldots,9).$ In this way, we have tested 1000 examples in one experiment and found that the King-like method (\ref{3}) fell in 643 examples. Clearly, the number of failures  depends on the choice of the set $\{a_0,a_1,a_2,a_3,a_4,a_5\}$ of random coefficients. For this reason we have performed 20 experiments with 1000 examples (in total, 20\,000 examples) and found that  the percentage of failures
ranged between 55\% and 65\%.

\medskip

Now we return to the group of mistakes (II). In essence, this  group of errors is irrelevant since the presented interval extensions do not guarantee the enclosure of the sought zero and we could omit any analysis of further mistakes (II).  However, we will discuss in short the convergence order given by theorems from the papers \cite{lotfi-1}--\cite{eft-numa},  as if the methods (\ref{3}) and (\ref{4}) had been properly constructed (see Remark 4). The following academic question arises: {\it If the methods {\rm (\ref{3})} and {\rm (\ref{4})} could preserve the inclusion property during the iterative process,
what would be the order of convergence of the semi-widths of produced inclusion intervals?}

\smallskip

Let $X=[\underline{x}, \overline{x}],$ and let $\rd(X)=(\overline{x}- \underline{x})/2$ be the semi-width of the interval $X.$ Similar denotation holds for the intervals appearing in   (\ref{3}) and (\ref{4}).
The determination of convergence order of the Ostrowski-like method (equal to 4) and the modified Ostrowski-like method (equal to 6) has been presented in \cite{ferara} and cited in \cite{eft-numa} for the same methods. These  iterative methods are given by (\ref{3}) and (\ref{4}) (respectively) for $\b=0.$ A very similar approach for finding the order of the King-like method (\ref{3}) ($\b\in \R$) (equal to 4) has been given in \cite{fazi} and cited in \cite{eft-AMC}. All of these results are wrong! The main reason for mistakes is
careless and incorrect application of properties of interval arithmetics, as mentioned at the beginning of this note. Although the mistakes of type (I), the lack of inclusion property, make  the ``interval" methods considered in  \cite{lotfi-1}--\cite{eft-numa} pointless, we give anyway the results of our numerical experiments concerned with the convergence order.

\smallskip

Take for the moment $c=1$ in the method (\ref{2}) realized in real arithmetic. Then the order of the iterative method is 3 (see \cite[Example 2.2]{elzever}). The use of corrective factor $c=(1+(\b-2) t)/(1+\b t)$ in (2) leads to the King method (\ref{2}) of order 4. Numerical experiments have confirmed that the corrective factor $c$ also accelerates the King-like method (\ref{3}) if this method runs correctly.
 However, the application of $c=(1+(\b-2) t)/(1+\b t)$ in King-like formula (\ref{3}) (defined in $I(\R)$) shifts the interval $F'(X)$ which can cause the loss of inclusion property.

 \smallskip

 To estimate the increase of the order of convergence of King-like method (\ref{3}) that produces inclusion intervals in each iteration, we have tested 20 functions for the parameter $\b=-2+0.5m\ (m=0,1,\ldots,9),$ furnishing in this way 200 examples. To measure convergence rate of the sequence of intervals $\{X_k\}$, we have applied  the so-called {\it computational order of convergence}  (see \cite{elzever})
 \be
 r_c=\frac{\log|\rd(X_{k+1})/\rd(X_k)|}{\log|\rd(X_k)/\rd(X_{k-1})|},\label{12}
 \ee
 which gives  very good results (quite close to the theoretical value of order of convergence)  if the radii involved in (\ref{12}) are very small.

 \smallskip

 To obtain considerably narrow intervals $X_k$ we have run 5 iterations and multi-precision arithmetic in computer algebra system {\it Mathematica} with 1000 significant digits, which simulates the limit process (appearing in finding theoretical convergence order) in a quite satisfactory manner. In the test examples when  the King-like method (\ref{3})  has worked well, we have been able to apply the formula (\ref{12}) and calculated $r_c=3.5+\eps,$ with $\eps=O(10^{-2})$ (disagreements on the second decimal digit). We have previously stressed that the determination of theoretical order of the ``interval" methods (\ref{3}) and (\ref{4}) is of negligible interest since these methods  are not defined well so that we have not analyzed the values of $r_c.$ In every case, the order cannot be four, as asserted in the mentioned papers.
 We end this  note with the observation that  numerical results presented in \cite{lotfi-1}, \cite{ferara} and \cite{eft-AMC} show that the order of convergence  is really slower than the order given by theorems in the cited papers, which means that theoretical results remain unconfirmed.

\bigskip
{\bf Acknowledgement.} This work was supported by the Serbian Ministry of Education and Science under Grant 174022.

\end{document}